\documentclass[11pt]{amsart}  \usepackage{amssymb,amsthm}

\newtheorem{theo}{Theorem}[section] \newtheorem{prop}[theo]{Proposition}
\newtheorem{coro}[theo]{Corollary} \newtheorem{lemma}[theo]{Lemma}
\newtheorem{teo}[theo]{Theorem} \newtheorem{lem}[theo]{Lemma}
\newtheorem{cor}[theo]{Corollary} \theoremstyle{definition}
\newtheorem{defi}[theo]{Definition} \newtheorem{rema}[theo]{Remark}
 
\newtheorem{remark}[theo]{Remark} 
 
\newcommand{\fin}{\hfill$\square$}

 \newcommand{\NN}{\mathbb{N}}
 
 \newcommand{\RR}{\mathbb{R}}
\renewcommand{\SS}{\mathbb{S}}

\newcommand{\la}{\langle}  \newcommand{\ra}{\rangle}
\newcommand{\tr}{\text{\rm tr}} \newcommand{\Ric}{\text{\rm Ric}}
\newcommand{\Riem}{\text{\rm Riem}} 
\newcommand{\II}{\text{\rm II}}

\title[Cosmological versus CMC time I: Flat spacetimes]{Cosmological time
  versus CMC time I: Flat spacetimes} 

\author[L. Andersson]{Lars Andersson$^\star$} \email{larsa@math.miami.edu}
\thanks{${}^\star$ Supported in part by the NSF, under  contract no. DMS
0104402 with the University of Miami.}  \address{Albert Einstein Institute,
Am M\"uhlenberg 1, D-14476 Potsdam, Germany \and Department of Mathematics,
University of Miami, Coral Gables, FL
 33124, USA}

\author[T. Barbot]{Thierry Barbot$^\dagger$}
\email{Thierry.BARBOT@umpa.ens-lyon.fr} 
\thanks{${}^\dagger\; {}^{\ddagger}   \; {}^\S$ 
Supported in
part by ACI ``Structures g\'eom\'etriques et Trous Noirs''.}

\address{CNRS, UMPA, \'Ecole Normale Sup\'erieure de Lyon}

\author[F. B\'eguin]{Fran\c cois B\'eguin$^\ddagger$}
\email{Francois.Beguin@math.u-psud.fr}  \address{Laboratoire de
Math\'ematiques, Universit\'e Paris Sud.}

\author[A. Zeghib]{Abdelghani Zeghib$^\S$}
\email{Abdelghani.ZEGHIB@umpa.ens-lyon.fr}  \address{CNRS, UMPA, \'Ecole
Normale Sup\'erieure de Lyon}

\date{April 28, 2006}

\begin{document}

\begin{abstract}
This paper gives a new proof that maximal, globally hyperbolic, flat
spacetimes of dimension $n\geq 3$ with compact Cauchy hypersurfaces are
globally foliated by Cauchy hypersurfaces of constant mean curvature,   and
that such spacetimes admit a globally defined constant mean curvature  time
function precisely when they are causally incomplete.  The proof, which is
based on using the level sets of the cosmological time function as barriers,
is conceptually simple and will provide the basis for  future work on
constant mean curvature  time functions in general constant curvature
spacetimes, as well for an analysis of the asymptotics of constant mean
foliations.
\end{abstract}

\maketitle



\section{Introduction} 
The study of the global properties of spacetimes solving the
Einstein equations
plays a central role in  differential geometry and  general
relativity. However, with the exception of results which rely on small data
assumptions (nonlinear stability results) or the assumption of symmetries,
many fundamental questions about the global structure of Einstein spacetimes
remain open, including cosmic censorship, structure of singularities, and
existence of global foliations by Cauchy hypersurfaces with controlled
geometry.  The Einstein equation is hyperbolic only in a weak sense, and
therefore in order to approach its Cauchy problem from a PDE point of view,
it is neccessary either to impose gauge conditions, or extract a hyperbolic
system by modifying the equation.  The constant mean curvature (CMC)
condition is an important gauge condition in the study of the Cauchy problem
of the Einstein equation, and hence in general relativity.  The CMC time
gauge is known to lead to a well-posed Cauchy problem in conjunction with the
zero shift condition \cite{ChB:ruggeri} as well as with the spatial harmonic
gauge condition \cite{AndMon}. In the Hamiltonian formulation of the Einstein
equation, the volume of a CMC hypersurface can be viewed as the canonical
dual to the CMC time, see \cite{Fischer-Moncrief}.  In the case of 2+1
dimensional spacetimes, this point of view leads to a formulation of the
Einstein equation in CMC gauge as  a time-dependent Hamiltonian system on the
cotangent bundle of Teichmuller space \cite{moncrief:teich}.

There are numerous results concerning the existence of global CMC foliations
and CMC time functions under various symmetry conditions, for spacetimes with
and without matter. See \cite{andersson:survey,rendall:survey} for recent
surveys.   It should be noted that examples of  Ricci flat spacetimes which
do not contain any CMC Cauchy hypersurface were recently constructed
\cite{chruscetal}. However, it is not yet known if these examples are stable.

Spacetimes with constant sectional curvature  constitute an important
subclass of spacetimes, where one may expect to understand the fundamental
questions,  including the cosmic censorship problem completely.  However,
even within this subclass, there are still open questions relating to  the
existence and properties of constant mean curvature foliations, and the
asymptotic structure at cosmological singularities  is not fully understood.

The systematic study of spacetimes of constant sectional curvature was
initiated by Mess \cite{mess}, following work by among others Margulis
\cite{margulis}  and Fried \cite{fried}.
%
The classification of maximal globally hyperbolic flat  spacetimes with
complete Cauchy hypersurfaces has recently been completed by Barbot
\cite{barflat}, following work of Bonsante  \cite{Bons2,bonsante} and others.

%

The purpose of this paper is to give a proof that maximal, globally
hyperbolic, flat spacetimes of dimension $n\geq 3$ with compact Cauchy
hypersurfaces are globally foliated by CMC Cauchy hypersurfaces,  and that
such spacetimes admit a global CMC time function precisely when they are
causally incomplete, see Theorem \ref{t.main-flat-case} below.  The proof is
based on using the level sets of the cosmological time function as barriers.
This result is not new, see remark \ref{rem:context}, but the method of proof
presented here is conceptually simple and will provide the basis for an
analysis CMC time functions in general constant curvature spacetimes, as well
as of the asymptotics of CMC foliations in future work.

Recall that a Lorentz manifold, or spacetime,  $(M,g)$ is \emph{globally
hyperbolic} if it contains a Cauchy hypersurface $S$, i.e. a weakly spacelike
hypersurface such that each inextendible Causal curve in $M$ intersects $S$.
The hypersurface $S$ may without loss of generality be assumed to be smooth
and strictly spacelike \cite{bernal:sanchez:2003,bernal:sanchez:2005}.
A globally hyperbolic spacetime is \emph{maximal} if it cannot be extended in
the class of globally hyperbolic spacetimes.  For brevity we use the acronym
MGHF for maximal, globally hyperbolic, flat spacetimes.
Let $S \subset M$ be a spacelike hypersurface in a spacetime of dimension $n$
and let $\nu$ be its future directed unitary normal. Then for $X,Y$ tangent
to $S$, the second fundamental form is given by  $\II(X,Y) = \la \nu,
\nabla_X Y\ra$.   The mean curvature of  $S$ is defined by $H = \tr\II/(n-1)$.
The hypersurface $S$ is said to have \emph{constant mean curvature} (CMC)  if
$H \big{|}_S$ is constant.  If $M$ satisfies the \emph{timelike convergence
condition} (i.e. if $\Ric(V,V) \geq 0$ for timelike vectors $V$) and has
compact Cauchy hypersurfaces, then for each $p \in M$ and for each $\tau \ne
0$,  there is at most one CMC surface containing $x$ with mean curvature
$\tau$.
A compact spacelike hypersurface in a globally hyperbolic spacetime is a
Cauchy hypersurface \cite{budic:etal}, so the leaves of a CMC foliation are
always Cauchy hypersurfaces if they are compact.
A time function $t: M \to I$ is a \emph{CMC time function} if the level sets
of $t$ are CMC Cauchy hypersurfaces  with $H(t^{-1}(\tau)) = \tau$ for all
$\tau \in I$.  In contrast to the situation for CMC hypersurfaces and
foliations,  a globally defined CMC time function with compact level sets is
unique, even if the timelike convergence condition fails to hold. The proof
is a straightforward application of the maximum principle, see \cite[\S
2]{BBZ} for details.

It is a basic fact that if an MGHF spacetime $(M,g)$ with compact Cauchy
hypersurfaces is causally complete, then  it is a quotient of the Minkowski
space $\RR^{1,n-1}$. In this case $M$ is foliated by flat, totally geodesic
Cauchy hypersurfaces. Therefore we may focus on the case when $M$ is causally
incomplete. Without loss of generality, assume that $M$ is past causally
incomplete.
Then $M$ is future complete,  and is the  quotient of a convex strict subset
of $\RR^{1,n-1}$  by a group of isometries acting freely and properly
discontinuously. This subset is in fact a future regular domain
$E^+(\Lambda)$, cf. definition \ref{def:regdom}. The past boundary, or Cauchy
horizon, of $E^+(\Lambda)$ represents in some sense the universal cover of
the past cosmological singularity of $M$.

The cosmological time function $\tau(p)$ is  defined as the maximal
Lorentzian  length of past directed causal curves starting at $p$. 
The cosmological time function 
is a $C^{1,1}_{\text{\rm Loc}}$ function, but not $C^2$ in general, 
and therefore the mean curvature of its
level sets must be interpreted in the weak sense, in terms of supporting
hypersurfaces.  
An analysis of the the weak mean curvature of the level sets
of the cosmological time function of flat spacetimes, 
and an application of the strong maximum principle of
\cite{gallomax}, enables us to show that the level sets of $\tau$ can be used
as barriers for CMC hypersurfaces. 
An important role in this analysis is played by the
notion of regular domain in $\RR^{1,n-1}$, introduced by Bonsante
\cite{bonsante}.

A future regular domain $E^+(\Lambda)$ is the intersection of the future of a
family $\Lambda$ of lightlike hyperplanes.  It can be shown that the
universal cover of a past causally incomplete MGHF spacetime is isometric to
a future regular domain. If $\Lambda$ has at least two elements, then
$E^+(\Lambda)$ has regular cosmological time function, in the sense that
$\tau$ is bounded from below and the limit of $\tau$ along past inextendible
causal geodesics is zero. In particular this is true for the universal cover
of an incomplete MGHF spacetime $M$, as well as for $M$ itself. See  \S
\ref{s.regular-domains} for details.  The level sets of $\tau$ have
interesting geometric properties.  Benedetti and Guadagnini
\cite{benedetti:guadagnini}  showed that in a 2+1 dimensional MGHF  spacetime
with compact Cauchy hypersurface of genus $> 1$, the geometry induced on the level
sets of $\tau$ precisely corresponds to a Thurston earthquake deformation
defined in terms of the holonomy data of $M$.

\subsection{Statement of results} 

We now state the main results in this paper. The first result characterizes
the generalized mean curvature of the level sets of the cosmological time
function in a regular domain.

\begin{teo}
\label{t.barriers}
Consider  a (future complete flat) regular domain $E^+(\Lambda)$ 
in $\RR^{1,n-1}$, and the
associated cosmological time $\tau: E^+(\Lambda) \rightarrow (0,
+\infty)$. Then, for every $a\in (0,+\infty)$, the level hypersurface $S_a =
\tau^{-1}(a)$ has generalized mean curvature bounded from below by 
$-\frac{1}{a}$, and from above by $-\frac{1}{(n-1)a}$.
\end{teo}

Our convention for second fundamental form and mean curvature are such that
the future hyperboloids in Minkowski space have negative mean curvature with
respect to the future directed normal, see section 
\ref{s.construction-barriers}. Clearly, Theorem \ref{t.barriers} holds for
quotients of regular domains, and such spaces therefore have crushing
singularity, since the level sets of the cosmological time function provide a
sequence of Cauchy hypersurfaces with uniformly diverging mean curvature.

For the case of spacetimes with compact Cauchy hypersurface,  a standard barrier
argument  yields existence of a CMC foliation.

\begin{theo}
\label{t.main-flat-case}
Let $(M,g)$ be a MGHF spacetime with compact
Cauchy hypersurfaces.
\begin{enumerate} 
\item If $(M,g)$ is both past and future geodesically complete then it does not
  admit any globally defined CMC time function, but it admits a unique CMC
  foliation.
\item \label{point:main} If $(M,g)$ is future geodesically complete, then it admits a globally defined
CMC time function  $\tau_{cmc}:M\rightarrow I$ where $I=(-\infty,0)$.
Furthermore, the  CMC and cosmological times are comparable: 
$$
\tau  \leq -\frac{1}{\tau_{cmc}} \leq (n-1)\tau.
$$
\item A similar statement, but with a time range $I=(0,+\infty)$, is true in  the
past geodesically  complete case.
\end{enumerate} 
In all cases, these foliations are analytic.
\end{theo}

\begin{remark}\label{rem:context} 
This result is not new. It was proved in \cite{AMT} in the 2+1 dimensionsal
case, assuming the existence of one CMC Cauchy hypersurface. In 
\cite{andflat}, a proof was given for 
the case of spacetimes with hyperbolic spatial topology.  
Finally, it has been observed in \cite{barflat}, that the general case
follows from the classification of MGHF spacetimes with compact Cauchy
hypersurfaces.  

The proof provided here is conceptually much simpler that the arguments given
in the above mentioned 
papers. More importantly, this proof can be adapted to the general constant
curvature case.  
The proof of the main part of Theorem \ref{t.main-flat-case}, the case when
$M$ is causally incomplete, makes use of the level sets of the cosmological
time function of the universal cover of $M$, which is a regular domain,  
as barriers in the construction of CMC hypersurfaces. In principle, this idea
generalizes immediately also to the case of constant non-zero curvature. 
However, the geometry
and global causality in the non-flat case are  sufficiently complicated that
the technical details require a separate paper \cite{ABBZ}.  There, we will
in particular investigate the structure of non-flat regular domains.

Further, the level sets of the cosmological
asymptotic behavior of the level sets of the cosmological time function is
intimately related to the geometry of the singularity itself, i.e. the
boundary of the universal cover of the spacetime. This will enable us in a 
forthcoming paper to analyze the asymptotic behavor of the CMC foliation at
the cosmological singularity of constant curvature spacetimes, see
\cite{ABBZ:asymptotic}. In particular, in the case of flat 
spacetimes, we are able to prove in \cite{ABBZ:asymptotic} the conjecture of
Benedetti and Guadagnini \cite{benedetti:guadagnini} that the limit of the
geometry of the level sets of the CMC time function in the Gromov sense is
the same as the limit of the geometry of the level sets of the cosmological
time function. In the 2+1 dimensional case, this limit can be identified as a
point on the Thurston boundary of Teichmuller space. While one expects the
limiting geometry of the cosmic time levels to be the same as the CMC time
levels in general, there is not yet a clear identification of the limiting
geometry except in the 2+1 dimensional flat case. 
\end{remark}

\begin{remark} 
There is no compactness condition on Cauchy hypersurfaces
in  Theorem  \ref{t.barriers}. However, a direct proof of existence  of CMC
hypersurfaces given barriers requires compactness. In a noncompact situation,
it is necessary to consider a sequence of Plateau problems, following ideas
developed in \cite{Treibergs}. It is natural to ask whether  any flat regular
domain has a CMC foliation.  In particular, given two level hypersurfaces of the
cosmological time function with mean curvatures bounded above and below by
$c$,  is there a CMC hypersurface with mean curvature $c$ between them?  
Similarly, given an
isometry group of a regular domain, does there exist CMC hypersurfaces, or
CMC foliations, invariant under the isometry group action? 
\end{remark}

\subsection*{Overview of the paper} 

The proof of
Theorem  \ref{t.barriers} is given in section \ref{s.construction-barriers},
which is the central part of this article.   
In the preceding sections, we review introduce some notions and preliminary
results which will be needed there. In section \ref{s.regular-domains}, some
basic facts about regular domains are recalled. The results here are mainly
due to Bonsante \cite{bonsante}. The definition and properties
of the cosmological time are given in section \ref{s.cosmological-time}.
The classification of
MGHF spacetimes with compact Cauchy hypersurface is given in section
\ref{s.Cauchy-compact}.  
Section \ref{s.singularity} discusses the past horizon, and the retraction to
the singularity of a future complete regular domain. 
In  \S \ref{foliations}, we will explain how to get from hypersurfaces with
prescribed mean curvature to a CMC foliation.  This technique is well known
to experts in the field, but since the details are somewhat scattered in the
literature,  we include them for the convenience of the reader. Along the
way, we also check that this works with our notion of generalized mean
curvature. In
particular, in the literature the strong energy condition is often assumed,
but we consider also the case of positive curvature (corresponding to
spacetimes of deSitter type), for future use in \cite{ABBZ}. Finally, in
section \ref{s.proof-main-flat-case} we give the proof of Theorem
\ref{t.main-flat-case}.

\section{Flat regular domains} 
\label{s.regular-domains}

Regular domains in Minkowski spacetime $\RR^{1,n-1}$ were first  defined by
F. Bonsante in~\cite {Bons2, bonsante} (generalizing a construction of
G. Mess in the  2+1-dimensional case, see \cite{mess}). Here we will use an
equivalent definition introduced in \cite{barflat}, since it appears to be
slightly more adapted to our purpose. For more details, we refer to section
$4.1$ of \cite{barflat}.

The importance of flat regular domains for our purpose comes from the fact
that they have regular cosmological time function, see Proposition
\ref{pro.regular}, and that each MGHF spacetime with compact, or more
generally complete, Cauchy hypersurface is a quotient of a flat regular domain,
see Theorem \ref{teo.dscompact}. Thus, the analysis of the singularity of
MGHF spacetimes can be carried out by studying the past boundary of flat
regular domains. This will be carried out in section \ref{s.singularity}.

\begin{defi}
The Penrose boundary ${\mathcal J}_{n-1}$ of the Minkowski spacetime
$\RR^{1,n-1}$ is the space of null affine hyperplanes of ${\mathbb R}^{1,n-1}$
\end{defi}

Let $N$ be an auxiliary euclidean metric on ${\mathbb R}^{1,n-1}$.  Let
${\mathcal S}^{n-2}$ be the set of future oriented null elements of ${\mathbb
R}^{1,n-1}$ with $N$-norm $1$. Then the map  which associates to a pair
$(u,a)$ the null hyperplane $H(u,a) = \{ x  | \la x , u \ra = a
\}$ is a bijection between ${\mathcal S}^{n-2} \times {\mathbb R}$ and
${\mathcal J}_{n-1}$.  It defines a topology on ${\mathcal J}_{n-1}$, which
coincides with the topology of ${\mathcal J}_{n-1}$ as a homogeneous space
under the action of the Poincar\'e group;  ${\mathcal J}_{n-1}$ is then
homeomorphic to ${\mathbb S}^{n-2} \times {\mathbb R}$.

For every element $p$ of ${\mathcal J}_{n-1}$, we denote by $E^+(p)$ the
future of $p$ in  ${\mathbb R}^{1,n-1}$, and by $E^-(p)$ the past of of
$p$. If  $p$ is the null hyperplane $H(u,a)$, then $E^+(p) = \{ x  |  \la
x , u \ra < a \}$ and $E^-(p) = \{ x  |  \la x , u \ra > a
\}$.  They are half-spaces, respectively future-complete and past-complete.
For every closed subset $\Lambda$ of ${\mathcal J}_{n-1}$, we define
$$E^\pm(\Lambda) = \bigcap_{p \in \Lambda} E^\pm(p).$$

\begin{defi} \label{def:regdom} 
A closed subset $\Lambda$ of ${\mathcal J}_{n-1}$ is said to be
\textit{future regular} (resp. \textit{past regular}) if it contains at least
two elements and if
$E^+(\Lambda)$ (resp. $E^-(\Lambda)$) is non-empty.

A \textit{future complete flat regular domain} is a domain of the form
$E^+(\Lambda)$ were $\Lambda$ is a future regular closed subset of ${\mathcal
J}_{n-1}$.  Similarly, a \textit{past complete flat regular domain} is a
domain of the form $E^-(\Lambda)$ were $\Lambda$ is a past regular closed
subset of ${\mathcal J}_{n-1}$.   A \textit{flat regular domain} is a future
complete regular domain or a past complete regular domain.
\end{defi}

See \S $4.2$ of \cite{barflat} where it is proved in particular that this
definition of flat regular domains coincides with Bonsante's definition.

\begin{remark}
A past regular closed set $\Lambda$ is not necessarily future regular.
Actually, a closed subset of ${\mathcal J}_{n-1}$ is past regular and future
regular if and only if it is compact (and contains at least two points). See
Corollary $4.11$ in \cite{barflat}.
\end{remark}

\begin{remark} 
In the rest of the paper, we will mainly be dealing with a
 past incomplete, future complete spacetimes, and many statements have an
 obvious time reversed analog. In the following we will not make any explicit
 statements concerning the time reversed situation, and leave it to the
 reader to rephrase the relevant definitions and results. 
\end{remark}

\subsection{Cosmological time}
\label{s.cosmological-time}

In any spacetime $(M,g)$, we can define the \textit{cosmological time\/}
(see \cite{cosmic}):

\begin{defi}
The cosmological time of a spacetime $(M,g)$ is the function
$\tau:M\rightarrow [0,+\infty]$ defined by
$$\tau(x)=\mbox{Sup}\{ L(\gamma) \mid \gamma \in {\mathcal R}^-(x) \},$$
where ${\mathcal R}^-(x)$ is the set of all past-oriented  causal curves
starting at $x$, and $L(\gamma)$ the lorentzian length of the causal curve
$\gamma$.
\end{defi}

In general, this function has a very bad behavior: for example, if $(M,g)$ is
Minkowski spacetime, then $\tau(x)=+\infty$ for every $x$.

\begin{defi}
\label{d.regular}
A spacetime $(M,g)$ is said to  have \textit{regular cosmological time\/} if
\begin{enumerate} 
\item $M$ has \textit{finite existence time,\/} i.e. $\tau(x) < +\infty$ for every  $x$ in $M$,
\item for every past-oriented inextendible causal curve $\gamma: [0, +\infty)
 \rightarrow M$, $\lim_{t \to \infty} \tau(\gamma(t)) = 0$.
\end{enumerate} 
\end{defi}

The following result gives a charaterization of spacetimes with regular
cosmological time.

\begin{teo}[\protect{\cite[Theorem 1.2]{cosmic}}]
\label{teo.cosmogood}
If $(M,g)$ has regular cosmological time, then:
\begin{enumerate}
\item $M$ is globally hyperbolic,
\item The cosmological time $\tau$ is a time function, i.e. $\tau$ is
continuous and is strictly increasing along future-oriented causal curves,
\item for each $x$ in $M$ there is a future-oriented timelike ray  $\gamma:
[0, \tau(x)] \rightarrow M$ realizing the distance from the "initial
singularity", that is, $\gamma$  is a unit speed geodesic which is maximal on
each segment and satisfies:
\[ \gamma(\tau(x))) = x \;\;\; \tau(\gamma(t)) = t \]
\item $\tau$ is locally Lipschitz, and admits first and second derivative
almost everywhere.
\end{enumerate}
\end{teo}

One of the cornerstones of Bonsante's work on flat regular domains is the
following proposition:

\begin{prop}
\label{pro.regular}
Future complete flat regular domains have regular cosmological time.
\end{prop}

\begin{proof}
See \cite[Proposition 4.3 and Corollary 4.4]{bonsante}.
\end{proof}

\subsection{Maximal globally hyperbolic flat spacetimes with compact Cauchy
hypersurfaces}
\label{s.Cauchy-compact}

\begin{prop}
\label{pro.quotients}
Let $E^+(\Lambda)\subset\RR^{1,n-1}$ be a future complete 
flat regular domain.  Let $\Gamma$
be a discrete torsion free group of isometries of Minkoswki spacetime
$\RR^{1,n-1}$ preserving $E^+(\Lambda)$. Then, the action of $\Gamma$ on
$E^+(\Lambda)$ is free  and properly  discontinuous, and the quotient space
$M_{\Lambda}^+(\Gamma) = \Gamma\backslash E^+(\Lambda)$ is a globally
hyperbolic   spacetime with regular cosmological time.
\end{prop}

\begin{proof}[Sketch of proof]
The proof that the action is free and properly discontinuous can be found in
\cite[Proposition 4.16]{barflat}.
The cosmological time $\tau$ is obviously $\Gamma$-invariant. Hence, it
induces a map $\hat{\tau}$ on the quotient $M_\Lambda^+(\Gamma)$. Since
inextendible causal curves in  $M_\Lambda^+(\Gamma)$ are projections of
causal curves in $E^+(\overline{\Lambda})$, the cosmological time on the
quotient $M_\Lambda^+(\Gamma)$ is the map $\hat{\tau}$. It follows easily
that $M_\Lambda^+(\Gamma)$ has regular cosmological time.
\end{proof}

Conversely:

\begin{teo}
\label{teo.dscompact}
Every MGHF
with compact Cauchy hypersurfaces
is the quotient of a flat regular domain or of the entire Minkowski space by
a torsion-free discrete subgroup of isometries.
More precisely, let $(M,g)$ be a $n$-dimensional  MGHF
spacetime with compact Cauchy hypersurfaces.
\begin{enumerate}
\item If $(M,g)$ is not past (resp. future) geodesically complete, then
$(M,g)$ is  the quotient of a future (resp. past) complete regular domain in
$\RR^{1,n-1}$ by a torsion-free discrete subgroup of
$\mbox{Isom}(\RR^{1,n-1})$.
\item If $(M,g)$ is geodesically complete then it is the quotient of
$\RR^{n-1,1}$ by a subgroup of $\mbox{Isom}(\RR^{1,n-1})$ containing a finite
index free abelian subgroup  generated by $n-1$ spacelike translations.
\end{enumerate}
\end{teo}

\begin{proof} 
It follows from the classification of MGHF spacetimes with compact Cauchy
hypersurfaces given in 
\cite{barflat}. The result in \cite{barflat} is more precise: it
characterizes up to finite index the possible torsion-free discrete subgroups.
\end{proof}

\begin{remark}
The natural setting for a result like Theorem~\ref{teo.dscompact} is not
really spacetimes with compact Cauchy hypersurfaces, but rather MGHF  spacetimes
with complete Cauchy hypersurfaces.  Indeed, every flat regular domain admits a
complete Cauchy hypersurface ( see \cite[Proposition 4.14]{barflat}).
Conversely, according to \cite[Theorem 1.1]{barflat}, every MGHF spacetime
with complete Cauchy hypersurface can be tamely embedded in the quotient of a flat
regular domain by a discrete group of isometries of Minkowski except if it is
geodesically complete or if it is \textit{an unipotent spacetime.\/}
Geodesically complete MGHF  spacetimes with complete Cauchy hypersurfaces are
quotients  of the entire Minkoswki space $\RR^{1,n-1}$ by a commutative
discrete group of spacelike translations. Flat unipotent spacetimes are
defined and described in \S $3.3$ of \cite{barflat} (see also \cite{fried});
every flat unipotent spacetime is the quotient of a domain
$\Omega\subset\RR^{1,n-1}$ by a unipotent discrete subgroup of
$\mbox{Isom}({\mathbb R}^{1,n-1})$,  where $\Omega$ is of one of the three
following forms: $\Omega=E^+(p)$, $\Omega=E^-(p)$ or $\Omega = E^{+}(p) \cap
E^-(p')$ where $p$ and $p'$ are two parallel null hyperplanes.
\end{remark}

\section{Past horizon and  initial singularity of a future complete flat regular domain}
\label{s.singularity}

In this section, we consider a future complete flat regular domain
$E^+(\Lambda)$. We will describe the past horizon, the initial singularity,
and the so-called "retraction to the initial singularity" of $E^+(\Lambda)$.

\subsection{Horizons}

According to Proposition~\ref{pro.regular} and  Theorem~\ref{teo.cosmogood},
$E^+(\Lambda)$ is globally hyperbolic.  Since $E^+(\Lambda)$ is a future
complete convex open domain in Minkowski space,  its boundary  ${\mathcal
H}^-(\Lambda)$ is a past horizon (and thus enjoys all the known  properties
of horizons).

Since ${\mathcal H}^-(\Lambda)$ is the boundary of a convex domain,  it
admits support hyperplanes at each of its points. And since $E^+(\Lambda)$ is
future complete,  the future in ${\mathbb R}^{1,n-1}$ of any point $p$ in
${\mathcal H}^-(\Lambda)$ is contained in $E^+(\Lambda)$.  But,   timelike
hyperplanes containing $p$ all intersect the future of  $p$, it then follows
that support hyperplanes  to ${\mathcal H}^-(\Lambda)$ are non-timelike.

\begin{lem}
\label{lem.flathori}
Let $p$ a point of the past horizon ${\mathcal H}^-(\Lambda)$ of a future
complete flat
regular domain $E^+(\Lambda)$. 
Let $C(p) \subset T_pX$ be the set of  future
oriented tangent vectors orthogonal to support hyperplanes to  ${\mathcal
H}^-(\Lambda)$ at $p$. Then $C(p)$ is the convex hull of its null elements.
Moreover, the null elements of $C(p)$ are precisely the normals to elements
of $\Lambda$ tangent to ${\mathcal H}^-(\Lambda)$ at $p$.
\end{lem}

\begin{proof}
See \cite[corollary 4.12]{bonsante} (see also \cite[Proposition 11]{mess}).
\end{proof}

\subsection{Retraction to the initial singularity}

According to point $(3)$ in Theorem~\ref{teo.cosmogood}, for every point $x$
in a flat regular domain there is a unique maximal timelike geodesic ray with
future endpoint $x$ realizing the "distance to the initial singularity": we
call such a geodesic ray a \textit{realizing\/} geodesic for $x$.

\begin{prop}
\label{pro.tout}
Let $x$ be an element of a future complete flat regular domain
$E^+(\Lambda)$. Then, there is an unique realizing geodesic for $x$.
\end{prop}

\begin{proof}
See \cite[Proposition 4.3]{bonsante}.
\end{proof}

\begin{defi}
\label{def.tight}
A unit speed future oriented timelike geodesic $\gamma: [0, T] \rightarrow
E^+(\Lambda)$ 
is tight if for every $t$ in $[0,T]$ the restriction $\gamma: [0,t]
\rightarrow E^+(\Lambda)$ is a realizing geodesic for $\gamma(t)$.
\end{defi}

\begin{prop}
\label{pro.tight}
Let $\gamma: [0, T] \rightarrow E^+(\Lambda)$ be an unit speed future oriented timelike
geodesic with initial point in the past horizon. Then the following
assertions are equivalent:
\begin{enumerate}
\item $\gamma$ is tight,
\item the derivative of $\gamma$ at $0$ is orthogonal to a  support
hyperplane at $\gamma(0)$.
\end{enumerate}
\end{prop}

\begin{proof}
See \cite[Proposition 4.3]{bonsante}.
\end{proof}

\begin{defi}
The initial singularity of a future complete flat regular domain
$E^+(\Lambda)$ is the set of points in the past horizon  admitting at least
two support hyperplanes; it will be denoted by
${\Sigma}^-(\Lambda)$. 
\end{defi}

\begin{prop}
\label{pro.rinsing}
The map which associates to any point $x$ of a regular domain $E^+(\Lambda)$,
the initial singularity of the unique realizing geodesic for $x$ is a
continuous map taking value  in $\Sigma^-(\Lambda)$. This map is denoted $r$,
and called ``retraction to the initial singularity".
\end{prop}

\begin{proof}
See \cite[Proposition 4.3 and 4.12]{bonsante}.
\end{proof}

\subsection{Description of the retraction map}

\begin{prop}
For every $p$ in the past singularity $\Sigma^-$, the preimage $r^{-1}(p)$ in
$E^+(\Lambda)$ is the union of complete timelike geodesic rays with initial
point at $p$.
\end{prop}

\begin{proof}
The Proposition is an immediate corollary of Proposition~\ref{pro.tout} and
\ref{pro.tight}.
\end{proof}

\begin{cor}
\label{cor.fiberopen}
Let $p$ be an element of the past horizon of $E^+(\Lambda)$ such that the
convex hull $C(p)$ of the null generators has non-empty interior in the space
of timelike tangent vectors at $p$. Then, $r^{-1}(p)$ is open in
$E^+(\Lambda)$.\fin
\end{cor}

\section{Cosmological levels as barriers, Proof of Theorem \ref{t.barriers}}
\label{s.construction-barriers}

If $S$ is a spacelike hypersurface in a spacetime $(M,g)$, then the
second fundamental form (also known as the extrinsic curvature) of 
$S$ at a point $x$ is defined as  
$\II(X,Y) = \la \nu , \nabla_{X}Y \ra = - \la \nabla_X\nu
, Y \ra$ where $X$, $Y$  are tangent vectors to $S$ at $x$ and $\nu$ is the
future oriented timelike normal of $S$  (with lorentzian norm $-1$).  The
mean curvature is defined in terms of the trace of $\II$ with respect to the
induced metric as $H_S = \tr \II/(n-1)$.  
This definition requires $S$ to be at least $C^2$. 
Nevertheless, in certain cases, one can give a meaning to the
assertion ``a topological hypersurface has mean curvature bounded from below
(or above) by some constant $c$''. A definition of this notion for rough
spacelike hypersurfaces was given in \cite[Definition 3.3]{gallomax}, making
use of the notion of supporting hypersurfaces with one-sided Hessian bound. 
The following
definition, which does not include the one-sided Hessian bound, 
is sufficient for our purposes in this paper. 
We will say that $S$ is a $C^0$-spacelike
hypersurface in $M$ if for each $x \in S$, there is a neighborhood $U$ of $x$
so that $S \cap U$ is edgeless and acausal in $U$, see \cite[Definition
  3.1]{gallomax}. 
  
\begin{defi}
\label{d.generalized-curvature}
Let $S$ be a $C^0$-spacelike hypersurface in a
spacetime $(M,g)$. Given a real number $c$, we
will say that $S$ \textit{has generalized mean curvature bounded from above
by $c$ at $x$},  denoted $H_S (x) \leq c$, if 
there is a geodesically convex
open neighborhood  $V$ of $x$ in $M$ and a smooth 
spacelike
hypersurface $\SS_x^-$ in $V$ such that~:
\begin{itemize}
\item[--] $x\in\SS_x^-$ and $\SS_x^-$ is contained in the past of  $S \cap V$
(in $V$),
\item[--] the mean curvature of $\SS_x^-$ at $x$ is bounded from above by $c$.
\end{itemize}

Similarly, we will say that $S$ \textit{has generalized mean curvature is
bounded from below by $c$ at $x$}, denoted $H_S(x) \geq c$, if,  
there is a geodesically convex open neighborhood $V$ of $x$ in $M$ and a smooth
spacelike hypersurface $\SS_x^+$ in $V$ such that~:
\begin{itemize}
\item[--] $x\in\SS_x^+$ and $\SS_x^+$ is contained in the past  of $S \cap V$
(with respect to   $V$),
\item[--] the mean curvature of $\SS_x^+$ at $x$ is bounded from below by $c$.
\end{itemize}
We will write $H_S \geq c$ and $H_S \leq c$ to denote that $S$ has
generalized mean curvature bounded from below respectively above by $c$ for
all $x \in S$. 
\end{defi}

\begin{rema}
Let $S$ be a smooth spacelike hypersurface in a spacetime $(M,g)$, and $c$ be
a real number. If $H_S \leq c$ or $H_S \geq c$ in the sense of the definition
above, then the maximum principle, see Proposition \ref{p.Harnack} below,
implies that the same bounds hold in terms of the usual sense. 
\end{rema}

\begin{rema}
Let $S$ be a $C^0$-spacelike hypersurface,  
and let $x$ be a point of $S$.  
Assume that there exists two numbers $c^-,c^+$ such that $S$
has generalized mean curvature bounded from below by $c^-$ and from above by
$c^+$ at $x$.  Then $S$ has a tangent plane at $x$. Indeed, the point $x$
belong to two smooth hypersurfaces $\SS_x^-$ and $\SS_x^+$ which are
(locally) respectively in the past and in the future $S$. In particular,
$\SS_x^-$ is locally in the past of $\SS_x^+$. This implies that the tangent
hyperplane of $\SS_x^-$ at $x$ coincides with the tangent hyperplane of
$\SS_x^+$. And since $S$ is between $\SS_x^-$ and $\SS_x^+$, this hyperplane
is also tangent to $S$.
\end{rema}

Let us recall the statement of Theorem \ref{t.barriers}:

\begin{teo}
\label{t.barriers-flat}
Consider  a future complete flat regular domain $E^+(\Lambda)$ and the
associated cosmological time $\tau: E^+(\Lambda) \rightarrow (0,
+\infty)$. Then, for every $a\in (0,+\infty)$, the hypersurface $S_a =
\tau^{-1}(a)$ has generalized mean curvature 
satisfying $-\frac{1}{a} \leq H_{S_a} \leq -\frac{1}{(n-1)a}$. \end{teo}

\begin{rema} 
What is  important in the proof of Theorem~\ref{t.main-flat-case} is just the fact
  that the hypersurface $S_a$ has generalized mean curvature 
satisfying $\alpha(a) \leq H_{S_a} \leq \beta(a)$, where
  $\alpha(a),\beta(a)\to -\infty$ when $a\to 0$, and $\alpha(a),\beta(a)\to
  0$ when $a\to +\infty$.
\end{rema}

\begin{proof}
Let $x$ be a point on the level set $S_a$. We denote by $\gamma: [0, a]
\rightarrow E^+(\Lambda)$ the unique realizing geodesic for $x$, with initial
point $p = r(x)$. Let $v$ be the future oriented unit speed tangent vector of
$\gamma$ at $p$. We denote as before by $C(p)$ the set of vectors in $T_p X$
orthogonal to support hyperplanes of the past horizon at $p$.

\subsection*{Construction of ${\mathbb S}^+_x$.}
Define ${\mathbb S}^+_x$ as the hyperboloid $\{ z | d(p,z) = a \}$.  Since
$E^+(\Lambda)$  
is geodesically convex, for any $z$ in ${\mathbb S}^+_x$   the
timelike geodesic $(p, z)$ is contained in $E^+(\Lambda)$. Hence, its length
$a$ is less than $\tau(z)$. The unique realizing geodesic for $z$ must
therefore intersect $S_a$. Hence, ${\mathbb S}^+_x$ is contained in the
future of $S_a$.  The tangent hyperplane to ${\mathbb S}^+_x$ at $x$ is the
hyperplane orthogonal to $c$ at $x$. Hence, ${\mathbb S}^+_x$ is tangent to
$S_a$ at $x$. Finally, the mean curvature of ${\mathbb S}^+_x$ is obviously
$-\frac{1}{a}$ everywhere. As a consequence, $S_a$ has generalized mean
curvature satisfying $H_{S_a} \geq -\frac{1}{a}$.

\subsection*{Construction of ${\mathbb S}^-_x$.}
According to Lemma~\ref{lem.flathori}, the tangent vector $v$ of the
realizing geodesic $\gamma$ introduced above, belongs to the convex hull
$C(p)$. Let $B$ be a finite subset of the null elements of $C(p)$ such that
$v$ lies in the convex hull of $B$. We choose moreover $B$ minimal, i.e.
such that for any proper subset $B' \subset B$, $v$ does not belong to the
convex hull of $B'$. An equivalent statement is that $v$ belongs to the
relative interior $\mbox{Conv}(B)$.

The null hyperplanes $p + w^\perp$ for $w$ in $B$ form a finite subset
$\Lambda_B$ of $\Lambda$.  Observe that since the convex hull of $B$ contains
the timelike vector $v$, $B$ contains at least two elements. Hence,
$E^+(\Lambda_B)$ is a future complete flat regular domain.

Obviously, $E^+(\Lambda_B)$ contains $E^+(\Lambda)$.  Hence ${\mathcal
H}^-(\Lambda_B)$ is contained in the causal past of $E^+(\Lambda)$. Moreover,
$E^+(\Lambda_B)$ contains the timelike geodesic  $\gamma$, and also $x$, and
its past horizon ${\mathcal H}^-(\Lambda_B)$ contains $p$. According to
Lemma~\ref{lem.flathori}, support hyperplanes to ${\mathcal H}^-(\Lambda_B)$
at $p$ are hyperplanes orthogonal to vectors in the convex hull of $B$.  In
particular, the hyperplane orthogonal to the timelike vector $v$ is a
spacelike  support hyperplane. It follows that $\gamma$  is a realizing
geodesic for $x$ in  $E^+(\Lambda_B)$. Hence, $\tau_B(x) = a$, where $\tau_B$
is the cosmological  time for $E^+(\Lambda_B)$.

Let $S'_B$ be the level set  $\{ \tau_B = a \}$ in $E^+(\Lambda_B)$, and
define ${\mathbb S}^-_x$ as a small open neighborhood  of $x$ in $S'_B \cap
E^+(\Lambda)$. Let $V$ be a geodesically convex neighborhood of $x$
containing ${\mathbb S}^-_x$ (for example, the Cauchy development of
${\mathbb S}^-_x$ in $E^+(\Lambda)$). For any $z$ in $S'_B$ let $c$ be the
unique realizing geodesic for $z$ in $E^+(\Lambda)$. Since ${\mathcal
H}^-(\Lambda_B)$   is in the causal past of ${\mathcal H}^-(\Lambda)$ there
is a past extension of $c$ with past endpoint in ${\mathcal H}^-(\Lambda_B)$.
Hence, $\tau(z) \leq a$. It follows that ${\mathbb S}^-_x$ lies in the causal
past of $S_a$ in $V$.

To complete the proof, we must prove  that ${\mathbb S}_x^-$ near $x$ is
smooth, admits at $x$ the same tangent hyperplane $(x-p)+v^\perp$, and that
it has constant mean curvature $-\frac{d}{(n-1)a}$ for some integer $1 \leq d
\leq n-1$.

Consider ${\mathbb R}^{1,n-1}$ as a vector space, with origin $p = 0$.  Let
$F$ be the vector space spanned by $\mbox{Conv}(B)$. Then $F$  is a timelike
subspace, with dimension $2 \leq k \leq n$,  and we have a splitting
${\mathbb R}^{1,n-1} = F \oplus F^\perp$.  The subspace $F^\perp$ is
spacelike. Every element of $\Lambda_B$ is a null hyperplane containing
$F^\perp$. It follows easily that  $E^+(\Lambda_B)$ is the sum $E'(\Lambda_B)
\oplus F^\perp$, where $E'(\Lambda_B) = F \cap E^+(\Lambda_B)$.  For every
element $H$ of $\Lambda_B$, $H \cap F$ is a null hyperplane in $F \approx
{\mathbb R}^{1,k-1}$.

Let $\Lambda'_B = \{ H \cap F \mid H \in \Lambda_B\}$.  Then $\Lambda'_B$ is
a finite subset of the Penrose boundary of $F$.  Clearly $E'(\Lambda_B)$ is
precisely the flat regular domain $E(\Lambda'_B) \subset F$. Now we observe
that restricting to $F$, $v$ is in the
interior of $\mbox{Conv}(\Lambda'_B)$. Hence, for some small neighborhood
$V'$ of $x$ in $E(\Lambda'_B)$, which can be selected  geodesically convex,
the image by the retraction $r$ of each point $y$ in $V'$  is $p$. Shrinking
$V$ if necessary, we can assume that $V$ is contained in  $V' \oplus
F^\perp$. According to Corollary~\ref{cor.fiberopen},   ${\mathbb S}^-_x$ has
the form ${\mathbb H} \oplus F^\perp$, where $\mathbb H$ is the hyperboloid
consisting of points in $F$ in the future of $p$ and at lorentzian distance
$a$ from $p$.  Hence, ${\mathbb S}^-_x$ is smooth and admits at $x$ the same
tangent hyperplane than $S_a$ (the orthogonal $x+v^\perp$). Moreover, since
$x+F^\perp$ is totally geodesic, and since  the principal directions of
$\mathbb H$ are all equal to $-\frac{1}{a}$, the mean curvature of ${\mathbb
S}^-_x$ is equal to $-\frac{1}{a}.\frac{d}{n-1}$ where $d=k-1$. This shows
that $S_a$ has generalized mean curvature satisfying $H_{S_a} \leq
-\frac{1}{(n-1)a}$. 
\end{proof} 

\begin{remark} \label{rem:strong}  
The proof of theorem \ref{t.barriers-flat} shows that the
  second fundamental forms of ${\mathbb S}^-_x, {\mathbb S}^+_x$ have
  eigenvalues $-1/a,0$ (in the case of ${\mathbb S}^-_x$) and $-1/a$ (in the
  case of ${\mathbb S}^+_x$. Therefore the level sets of $\tau$ have mean
  curvature satisfying $-1/a \leq H_{S_a} \leq -1/((n-1)a)$ with one-sided 
Hessian
  bound as in \cite[Definition 3.3]{gallomax}, and hence the strong maximum
  principle for spacelike hypersurfaces given in \cite[Theorem
  3.6]{gallomax} applies in our situation. However, we shall not need
  the full strength of this result here. See proposition \ref{p.Harnack}
  below for the version of the maximum principle which we shall make use of. 
\end{remark} 

The eigenvalue bounds stated in remark \ref{rem:strong} allow us to
give a more precise characterization of  the regularity of the cosmological
time function. The Hessian bounds for the height function implied by the
bounds on the second fundamental form of the supporting hypersurfaces,
together with an application of the case $p=\infty$ of 
\cite[Proposition 1.1]{caffarelli:cabre}  proves
\begin{coro} \label{cor:C11} 
$\tau \in C^{1,1}_{\text{\rm Loc}}$ 
\end{coro}

We leave it to the reader to formulate the obvious analogs of 
theorem \ref{t.barriers-flat} and corollary
\ref{cor:C11} for \emph{past complete} flat regular domains $E^-(\Lambda)$ 
which hold in terms of the 
\emph{reverse} cosmological time $\widehat{\tau}: E^-(\Lambda) \to
(0,+\infty)$. 

%




\section{From barriers to CMC time functions}
\label{foliations}

In this section, we consider a $n$-dimensional, $n \geq 3$, maximal globally
hyperbolic spacetime $(M,g)$ with compact Cauchy hypersurfaces and 
constant curvature equal to $k$. We emphasize that many of the proofs that we
give are not valid without the assumption that $M$ has compact Cauchy
surfaces. Recall
that $(M,g)$ has curvature $k$ if the Riemann tensor satisfies 
$$
\la \Riem(X,Y)Y,X \ra = k ( \la X, X \ra \la Y , Y \ra - \la X , Y \ra^2 ) 
$$
for any vector fields $X,Y$. Then the Ricci tensor satisfies 
$\Ric = (n-1) k g$.  
We will define a notion of \textit{sequence of asymptotic barriers}, and
prove (using  quite classical arguments) that $(M,g)$ admits a CMC time function provided  that it admits a sequence of asymptotic barriers.

\begin{defi}
Let $c$ be a real number. A \textit{pair of $c$-barriers} is a
pair of $C^0$-spacelike Cauchy hypersurfaces $(\Sigma^-,\Sigma^+)$ in $M$ such that
\begin{itemize} 
\item[--] $\Sigma^+$ is in the future of $\Sigma^-$,
\item[--] $H_{\Sigma^-} \leq c \leq H_{\Sigma^+}$
in the sense of definition \ref{d.generalized-curvature}.
\end{itemize} 
\end{defi}

\begin{defi}
Let $\alpha$ be a real number.  A \textit{sequence of asymptotic past
$\alpha$-barriers} is a sequence of $C^0$-spacelike Cauchy hypersurfaces
$(\Sigma_m^-)_{m\in\NN}$ in $M$ such that
\begin{itemize}
\item[--] $\Sigma_m^-$ tends to the past end of $M$ when $m\to
  +\infty$ (\textit{i.e.} given any compact subset $K$ of $M$, there exists $m_0$ such that
  $K$ is in the future of $\Sigma_m^-$ for every $m\geq m_0$), 
\item[--] $a_m^- \leq H_{\Sigma_m^-} \leq a_m^+$, 
where $a_m^-$ and $a_m^+$ are real numbers
  such that $\alpha<a_m^-\leq a_m^+$, and such that $a_m^+\rightarrow \alpha$
 when $m\to +\infty$.  
\end{itemize}

Similarly, a \textit{sequence of asymptotic future $\beta$-barriers} is a
sequence of $C^0$-spacelike Cauchy hypersurfaces $(\Sigma_m^+)_{m\in\NN}$
in $M$ such that
\begin{itemize}
\item[--] $\Sigma_m^+$ tends to the future end of $M$ when $m\to +\infty$,
\item[--] $b_m^- \leq H_{\Sigma_m^+} \leq b_m^+$, 
where $b_m^-$ and $b_m^+$ are real numbers
such that $b_m^-\leq b_m^+<b$, and such that $b_m^-\rightarrow \beta$
 when $m\to +\infty$. 
\end{itemize}
\end{defi}

\begin{theo}
\label{t.foliation} Let $(M,g)$ be an 
$n$-dimensional, $n \geq 3$, maximal globally
hyperbolic spacetime, with compact Cauchy hypersurfaces and 
constant curvature $k$, and such that 
$(M,g)$ admits a sequence of  asymptotic past
$\alpha$-barriers and a sequence of asymptotic future $\beta$-barriers.
If $k\geq 0$, assume moreover that $(\alpha,\beta)\cap [-\sqrt{k},\sqrt{k}]=\emptyset$. Then, $(M,g)$ admits a CMC time function $\tau_{cmc}:M\rightarrow (\alpha,\beta)$.   
\end{theo}

Theorem \ref{t.foliation} follows easily from known facts in case the
barriers are smooth, and introducing $C^0$ barriers is not difficult given the
results above. Nevertheless, 
since we are not aware of a reference for this precise statment, we include a
proof below. The following are the two main technical steps in the proof. In
the case of smooth barriers and hypersurfaces, they were proved in this 
formulation by Gerhardt \cite{Ger3}. 

\begin{itemize}
\item[--] a proposition which states that any CMC
  hypersurface of mean curvature $c'$ lies in the future of any CMC
  hypersurface of mean curvature $c$ whenever $c'>c$
  (Proposition~\ref{p.order});   
\item[--] a theorem 
which ensures the
  existence of a Cauchy hypersurface of constant mean curvature $c$,
  assuming the existence of  a pair of $c$-barriers
  (Theorem~\ref{t.CMC-existence}). 
\end{itemize}
Let us start with a slight generalization of the classical \textit{maximum
  principle}.  

\begin{prop}
\label{p.Harnack}
Let $\Sigma$ and $\Sigma'$ be two $C^0$-spacelike hypersurfaces. 
Assume that these hypersurface have one point
$x$ in common, and assume that $\Sigma$ is in the past of $\Sigma'$. 
Assume that $\Sigma$ has
generalized mean curvature bounded from above by $c$ at $x$, and $\Sigma'$ has
generalized mean curvature bounded from below by $c'$ at $x$. 
Then $c\geq c'$. 
\end{prop}

\begin{rema} 
Proposition \ref{p.Harnack}, which may be viewed as a
comparison principle, follows from the strong
maximum principle for $C^0$ hypersurfaces satisfying a one-sided Hessian
bound, see \cite[Theorem 3.6]{gallomax}. The notion of generalized mean
curvature we are using here does not included this requirement and we
therefore include the simple proof of the proposition. 
\end{rema} 

\begin{proof}
Since $\Sigma$ has generalized mean curvature bounded from above by $c$ at $x$,
there exists a smooth spacelike hypersurface $S_x$ such that $x\in S_x$,
$S_x$ is in the past of $\Sigma$ and the mean curvature of $S_x$ at
$x$ is at most $c$. Similarly, there exists a smooth spacelike
hypersurface $S_x'$ such that $x\in S_x'$, $S_x'$ is in the future of
$\Sigma'$ and the mean curvature of $S_x'$ at $x$ is at least
$c'$. Since $\Sigma$ is in the past of $\Sigma'$, this implies that
$S_x$ is in the past of $S_x'$. And since the point $x$ belongs to
both $S_x$ and $S_x'$, we deduce that $S_x$ and $S_x'$ share the same
tangent hyperplane at $x$. Now the classical maximum principle can be applied
to show that $c \geq c'$. 
\end{proof}

The following result was proved by Gerhardt for the case of spacetimes with
a lower bound on the Ricci curvature on timelike vectors, see \cite[Lemma
2.1]{Ger3}. 

\begin{prop}
\label{p.order}
 Let $(M,g)$ be an 
$n$-dimensional, $n \geq 3$, maximal globally
hyperbolic spacetime, with compact Cauchy hypersurfaces and 
constant curvature $k$.
Let $\Sigma$ and $\Sigma'$ be two smooth Cauchy hypersurfaces in
$M$. Assume that $H_\Sigma \leq c$ and $H_{\Sigma'} \geq c'$, with $c \leq
c'$.  
If $k$ is non-negative, assume moreover that
$c<- \sqrt{k}$ 
or that $c'>\sqrt{k}$. Then $\Sigma'$ is in the future of 
$\Sigma$.  
\end{prop}

We will give a proof of Proposition \ref{p.order} below, as we shall make use
of some of the details in the proof of theorem \ref{t.foliation}.

Let $\Sigma_0$ be a smooth
Cauchy hypersurface with future unit normal $\nu_0$. 
Recall that the orbit of the {\em Gauss flow} of smooth Cauchy hypersurface 
$\Sigma_0$ in the direction
$\nu_0$ consists of the Cauchy hypersurfaces $\Sigma_t = F_t(\Sigma_0)$ where $F:
I \times \Sigma_0 \to M$ is defined as $F_t(x) = exp_x(t \nu_0)$ for $x \in
\Sigma_0$, 
for $t \in I$. Here $I$ is the maximal time
interval where $F_t$ is regular. 
The core of the proof of Proposition~\ref{p.order} is the following standard
comparison lemma, see for example \cite[corollary 2.4]{AndHow}.

\begin{lemma}
\label{l.Gauss-flow}
We consider the orbit $(\Sigma_t)_{t\in I}$ of a smooth Cauchy
hypersurface $\Sigma_0$ under the Gauss flow. We consider a geodesic
$\gamma$ which is orthogonal to the $\Sigma_t$'s, and we denote by
$p(t)$ the point of intersection of the geodesic $\gamma$ with the
hypersurface $\Sigma_t$. 
The mean curvature $H(t)$ of $\Sigma_t$ at $p(t)$ satisfies the
differential inequality
$$\frac{dH(t)}{dt}\geq (n-1)( H(t)^2 - k).$$
\end{lemma}

\begin{proof}[Proof of Proposition~\ref{p.order}]
Assume that $\Sigma'$ is not in the future of $\Sigma$. Then, we can
consider a future-directed timelike geodesic segment $\gamma$ going
from a point of $\Sigma'$ to a point of $\Sigma$
having maximal length among all such geodesic segments. 
It is well-known that $\gamma$ is
orthogonal to both $\Sigma'$ and $\Sigma$, and that there is no focal
point to $\Sigma'$ or $\Sigma$ along $\gamma$ (see
e.g.~\cite[Proposition 4.5.9]{HawEll}). We will denote by $p'\in\Sigma'$ and $p\in\Sigma$ the 
ends of $\gamma$, and  by $\delta$ be the length of $\gamma$.

If $k$ is non-negative, we have to distinguish two differents cases,
according to whether $c'>\sqrt{k}$ or $c<-\sqrt{k}$. 
Let us consider the first case. Since there is no focal point to
$\Sigma'$ along $\gamma$, the image $\Sigma_t'$ of $\Sigma'$ by the
time $t$ of the Gauss flow is well-defined for $t\in [0,\delta]$ in a
neighbourhood of $\gamma$. Denote by $p'(t)$ the point of intersection
of the hypersurface $\Sigma_t'$ with the geodesic segment $\gamma$,
and by $H'(t)$ the mean curvature of $\Sigma_t'$ at $p'(t)$. By
Lemma~\ref{l.Gauss-flow}, $t\mapsto H'(t)$ satisfies the differential
inequality $\frac{dH'}{dt}\geq (n-1)(H'^2 - k).$ 
This implies that $H'$ increases along $\gamma$ (note that $H'(t)^2$
is strictly greater 
than $k$ for every $t$, since $H'(0)=c'>\sqrt{k}$ by assumption 
and since $H'(t)$ increases).  
In particular, we have $H'(\delta)>H'(0)=c'$. But now, recall
that, by definition of $\Sigma_\delta'$, every point of
$\Sigma_\delta'$ in a neighbourhood of $\gamma(\delta)=p$ is at distance
exactly $\delta$ of $\Sigma'$. Also recall that $\gamma$ is the
longest geodesic segment joining a point of $\Sigma'$ to a point of
$\Sigma$. This implies that $\Sigma$ is in the past of
$\Sigma_\delta'$. Hence, by Proposition~\ref{p.Harnack}, the mean
curvature of $\Sigma$ at $p$ is bounded from below by the mean
curvature of $\Sigma_\delta'$, which itself is strictly greater than the mean
curvature of $\Sigma'$. This contradicts the assumption $c\leq c'$. 

The proof is the same in the case where $c<-\sqrt{k}$ (except that one considers the backward orbit of $\Sigma$ for the Gauss flow, instead of the forward orbit of $\Sigma'$). 
\end{proof}

\begin{rema}
\label{r.uniqueness}
Proposition~\ref{p.order} implies that, for every $c\in\RR\setminus
[-\sqrt{k},\sqrt{k}]$, there exists at most one Cauchy hypersurface in
$M$ with constant mean curvature equal to $c$.   
In particular, for any open interval $(\alpha,\beta)$, which if $k\geq 0$
satisfies the condition $(\alpha,\beta)\cap [-\sqrt{k},\sqrt{k}]=\emptyset$,
there exists at most one function $t_{cmc}:M\to (\alpha,\beta)$ 
such that $t_{cmc}^{-1}(c)$ is a smooth
Cauchy hypersurface with constant mean curvature equal to $c$ for every $c\in
(\alpha,\beta)$. Note that we are not assuming here that $t_{cmc}$ is a time function  
(recall that, if $t_{cmc}$ is a time function, then it is  
automatically unique, without any assumption on $(\alpha,\beta)$).

Further, it is easy to see using a
maximum principle argument, that in the standard deSitter space with topology
$S^{n-1} \times \RR$ and curvature $k > 0$, there is no Cauchy hypersurface with
mean curvature $c \in \RR \setminus [-\sqrt{k},\sqrt{k}]$. Therefore
Proposition \ref{p.order} is vacuous in this case. 
\end{rema}

\begin{theo}
\label{t.CMC-existence} Let $(M,g)$ be an 
$n$-dimensional, $n \geq 3$, maximal globally
hyperbolic spacetime, with compact Cauchy hypersurfaces.
Let $c$ be any real number, and assume that there exists a pair of
$c$-barriers $(\Sigma^-,\Sigma^+)$ in $M$. Then, there exists a smooth
Cauchy hypersurface $\Sigma$ with
constant mean curvature equal to $c$. Moreover, $\Sigma$ is in the future of
$\Sigma^-$ and in the past of $\Sigma^+$.
\end{theo}

\begin{proof}
The result is proved e.g. in~\cite{Ger1} in the case where the barriers
$\Sigma^-$ and $\Sigma^+$ are smooth. The only way the
barriers $\Sigma^-$ and $\Sigma^+$ are used in Gerhardt's proof is via the
maximum principle (to show that a family of Cauchy hypersurfaces whose mean
curvature approaches $c$ cannot  ``escape to infinity''). 
Since the  maximum principle is still valid for $C^0$ hypersurfaces
(Proposition~\ref{p.Harnack}), Gerhardt's proof also applies in the case
where the barriers are not smooth.
\end{proof}

\begin{proof}[Proof of Theorem~\ref{t.foliation}]
We consider a sequence $(\Sigma_m^-)_{m\in\NN}$ of asymptotic past
$\alpha$-barriers, and a sequence $(\Sigma_m^+)_{m\in\NN}$ of
asymptotic future $\beta$-barriers.  

\subsection*{Construction of the function $\tau_{cmc}$.}
Fix $c\in (\alpha,\beta)$. For $m$ large enough, the
pair of Cauchy hypersurfaces $(\Sigma_m^-,\Sigma_m^+)$ is a pair of
$c$-barriers. Thus, by Theorem~\ref{t.CMC-existence}, for any $c \in
(\alpha,\beta)$, there exists a
Cauchy hypersurface $S_c$ with constant mean curvature equal to
$c$. Proposition~\ref{p.order} implies that the $S_c$'s are pairwise disjoint, 
and that $S_c$ is in the past of $S_{c'}$ if $c<c'$ (let us call this 
``property $(\star)$").

Now, let us prove that the set $\bigcup_{c\in (\alpha,\beta)} S_c$ is
connected. Assume the contrary. Because of property $(\star)$, there are only two possible cases~:  
\begin{itemize}
\item[(i)] there exists $c_0\in (\alpha,\beta)$ such that 
$\bigcup_{c>c_0} I^+(S_c)\subsetneq I^+(S_{c_0}),$
\item[(ii)] or there exists $c_0\in (\alpha,\beta)$ such that $\bigcup_{c<c_0} I^-(S_c)\subsetneq I^-(S_{c_0}).$
\end{itemize}
Let us consider, for example, case~(i). Using the Gauss flow, we
can push the hypersurface $S_{c_0}$ towards the future, in order to
obtain a Cauchy hypersurface $S_{c_0}'$ which is in the future of
$\Sigma_{c_0}$, but as close to $S_{c_0}$ as we want. In
particular, we can assume that $S_{c_0}'$ is not in the future of
$S_c$ for any $c>c_0$. Moreover, according to
Lemma~\ref{l.Gauss-flow}, the mean curvature of $S_{c_0}'$ is bounded
from below by some number $c_0'>c_0$. But this contradicts
Proposition~\ref{p.order}. Case~(ii) can be treated similarly. As a
consequence, the set $\bigcup_{c\in (\alpha,\beta)} S_c$ is
connected. Note that this implies that the hypersurface $S_c$ depends
continuously on $c$

Now, let us prove that the union $\bigcup_{c\in (\alpha,\beta)} S_c$
is equal to the whole $M$. Assume that there exists a point
$x\in M\setminus \bigcup_{c\in (\alpha,\beta)} S_c$. Since
the hypersurface $S_c$ depends continuously on $c$, there are
only two possible cases~:
\begin{itemize}
\item[(i)] either $x$ is in the future of $S_c$ for every $c\in
(\alpha,\beta)$, 
\item[(ii)] or $x$ is in the past of $S_c$ for every $c\in
(\alpha,\beta)$. 
\end{itemize}
Now, recall that we have a sequence $(S_m^+)_{m\in\NN}$ of
asymptotic future $\beta$-barriers. By definition, this means that
$S_m^+$ has generalized mean curvature bounded from below by some $b_m^-$ and
smaller than some $b_m^+$ where $b_m^-\leq b_m^+<\beta$ and 
$b_m^-\to_{m\to\infty}\beta$. Fix some integer $p$. One can find
$q>p$ such that $b_q^->b_p^+$. Then $(S_p^+,S_q^-)$ is a
pair of $b_p^+$-barriers. By Theorem~\ref{t.CMC-existence}, one can find a
Cauchy hypersurface with constant mean curvature equal to $b_p^+$
between  $S_p^+$ and $S_q^-$, and by uniqueness (see
remark~\ref{r.uniqueness}), this hypersurface is the hypersurface
$S_c$ for 
$c=b_p^+$. In particular, for $c\geq b_p^+$, the hypersurface
$S_c$ is in the future of the barrier $S_p^+$. Now, recall
that, by definition of a sequence of asymptotic future barriers,
$S_p^+$ tends to the future end of $M$ when $p\to\infty$. This
shows that case (i) cannot happen. Of course, one can exclude
case (ii) using similar arguments. Therefore we have proved that
$\bigcup_{c\in (\alpha,\beta)} S_c = M.$

Now, we can define define the function $\tau_{cmc}:m\to(\alpha,\beta)$
as follows~: for every $x\in M$, we set $\tau_{cmc}(x)=c$ where $c$ is
the unique number such that $x\in S_c$.

\subsection*{Properties of the function $\tau_{cmc}$.} The fact the
hypersurface $S_c$ depends continuously on $c$ implies that the
function $\tau_{cmc}$ is continuous. The fact that the hypersurface
$S_{c'}$ is in the strict future of the hypersurface $S_{c}$ when
$c'>c$ implies that the function $\tau_{cmc}$ is strictly increasing
along any future directed timelike curve. Hence,
$\tau_{cmc}$ is a time function.
\end{proof}

\begin{rema}
The function $\tau_{cmc}$ is also a time function in the following
stronger sense~: for every future directed timelike curve
$\gamma:I\rightarrow\RR$, one has
$$\frac{d}{dt}\tau_{cmc}(\gamma(t))>0.$$ 
Indeed, fix such a curve $\gamma$ and some $t_0\in I$, let
$x_0=\gamma(t_0)$ 
and $c_0=\tau_{cmc}(x_0)$. For $t$ small enough, denote by
$S_{c_0}^t$ the image of the hypersurface $S_{c_0}$ by the
time $t$ of the Gauss flow. Since the derivative $\gamma$ is
future-oriented timelike vector, there exists a constant $\lambda_1>0$ such
that, for $h>0$ small enough, the point $\gamma(t_0+h)$ is in the
future of the image of the hypersurface
$S_{c_0}^{\lambda_1.h}$. Now Lemma~\ref{l.Gauss-flow} implies that there
exists a constant $\lambda_2>0$ such that the mean curvature of the
hypersurface $S_{c_0}^{\lambda_1.h}$ is bounded from below by $c_0+\lambda_1.\lambda_2.h$ 
(for $h$ small enough). Then Proposition~\ref{p.order} implies that 
$S_{c_0+\lambda_1.\lambda_2.h}$ is in the past of $S_{c_0}^{\lambda_1.h}$. 
In particular, for $h$ small enough, the point $\gamma(t_0+h)$ is the future of the
hypersurface $S_{c_0+\lambda_1.\lambda_2.h}$. In other words, we
have $\tau_{cmc}(t_0+h)>c_0+\lambda_1.\lambda_2.h$. This implies 
$\frac{d}{dt}\tau_{cmc}(\gamma(t))>\lambda_1.\lambda_2>0$. 
\end{rema}

\begin{rema}
Using the same arguments as above, one can prove the following result:

\medskip

\noindent \textit{Let $(M,g)$ be an 
$n$-dimensional, $n \geq 3$, maximal globally
hyperbolic spacetime, with compact Cauchy hypersurfaces and 
constant curvature $k$. 
Assume that $(M,g)$  admits a sequence of asymptotic past
$\alpha$-barriers. If $k\geq 0$, assume moreover that $\alpha\notin [-\sqrt{k},\sqrt{k}]$. 
Then, $(M,g)$ admits a partially defined CMC time function $\tau_{cmc}:U\rightarrow (\alpha,\beta)$  
where $U$ is a neighbourhood of the past end of $M$ (i.e. the past of a Cauchy hypersurface 
in $M$) and $\beta$ is a real number greater than $\alpha$.}
\end{rema}

\begin{prop}
\label{p.analyticity}
Let $(M,g)$ be an 
$n$-dimensional, $n \geq 3$, maximal globally
hyperbolic spacetime, with compact Cauchy hypersurfaces and 
constant curvature $k$.
%
Suppose that there
exists a function $t_{cmc}:M\to (\alpha,\beta)$ such that $t_{cmc}^{-1}(c)$
is a Cauchy hypersurface with constant mean curvature equal to $c$ for every
$c\in (\alpha,\beta)$. Assume moreover that one of the following hypotheses
is satisfied: 
\begin{itemize}
\item $t_{cmc}$ is a time function,
\item the curvature $k$ is negative,
\item the curvature $k$ is non-negative and $(\alpha,\beta)\cap [-\sqrt{k},\sqrt{k}]=\emptyset$. 
\end{itemize}
Then $t_{cmc}$ is real analytic. 
\end{prop} 

\begin{proof}[Sketch of proof]
Under the stated conditions, there is exactly one CMC Cauchy hypersurface for each
$c \in (\alpha,\beta)$. CMC hypersurfaces in a real analytic spacetime are real
analytic, since they are solutions of a quasi-linear elliptic PDE. Given a
CMC Cauchy hypersurface $S_0$ with mean curvature $c_0 \in (\alpha,\beta)$, let
$u$ be the Lorentz distance to $S_0$. For 
$c$ close to $c_0$, a Cauchy hypersurface $S_c$ with mean curvature $c$ is a graph
over $S_0$, defined by the level 
function $w = u \big{|}_{S_c}$. The function $w$
solves the mean curvature equation $H[w] = c$, 
which is a quasilinear elliptic system with real 
analytic dependence on the coefficients. It follows that $S_c$ depends in a
real-analytic manner on $c$, and that 
the 
function $t_{cmc}$ is a real analytic function on $M$. 
\end{proof}

\section{Proof of  Theorem~\ref{t.main-flat-case}}
\label{s.proof-main-flat-case} 

Let $(M,g)$ be a $n$-dimensional MGHF 
spacetime with compact Cauchy hypersurface. 
We 
first consider the case where $(M,g)$ is not past geodesically
complete. Then Theorem~\ref{teo.dscompact} states that $(M,g)$ is the
quotient of a future complete regular domain
$E^+(\Lambda)\subset\RR^{1,n-1}$  by a torsion-free discrete subgroup
$\Gamma$ of $\mbox{Isom}(\RR^{1,n-1})$.  Let
$\tau:E^+(\Lambda)\rightarrow (0,+\infty)$ be the cosmological time of
$E^+(\Lambda)$. It follows from Theorem~\ref{teo.cosmogood} and its proof,
see \cite[Proposition 2.2]{cosmic}, that for every
$a\in (0,+\infty)$, the level set $S_a=\tau^{-1}(a)$ is 
a closed strictly achronal edgeless hypersurface in $E^+(\Lambda)$. 
Moreover,
$\tau$ is obviously invariant  under every element of
$\mbox{Isom}(\RR^{1,n-1})$ preserving $E^+(\Lambda)$. Hence, for every
$a\in (0,+\infty)$, the projection $\Sigma_a$ of $S_a$ in 
$M\equiv\Gamma\setminus E^+(\Lambda)$ is a closed strictly achronal
edgeless hypersurface 
in $M$.  Since $M$ is globally hyperbolic with compact 
Cauchy hypersurfaces, this
implies that $\Sigma_a$ is a compact strictly achronal hypersurface in
$M$, and thus is a topological Cauchy hypersurface in $M$. 
Theorem~\ref{t.barriers-flat} implies that, for every $a\in
(0,+\infty)$, $\Sigma_a$ has generalized mean curvature bounded from below
by
$-1/a$, and bounded from above by $-1/((n-1)a)$.
Let $(a_m)_{m\in\NN}$ be a decreasing sequence of positive real
numbers such that $a_m\to 0$ when $m\to +\infty$, and
$(b_m)_{m\in\NN}$ be a increasing sequence of positive real numbers
such that $b_m\to +\infty$ when $m\to +\infty$. Observe that
$(\Sigma_{a_m})_{m\in\NN}$ is a sequence of past asymptotic 
$\alpha$-barrier in $M$ for $\alpha=-\infty$ (indeed
$-\infty<-1/a_m<-1/((n-1)a_m)$ for every $m$, and since
$-1/((n-1)a_m)\to -\infty$ when $m\to\infty$),  and
$(\Sigma_{b_m})_{m\in\NN}$ is a 
sequence of future asymptotic $\beta$-barrier in $M$ for $\beta=0$
(indeed $-1/b_m<-1/((n-1)b_m)<0)$ for every $m$, and since $-1/b_m\to
0$ when $m\to\infty$). Hence Theorem~\ref{t.foliation} implies 
that $M$ admits a globally defined CMC time function $\tau_{cmc}:M\to
(-\infty,0)$.  

Next, we prove that $\tau$ and $\tau_{cmc}$ are comparable. 
It follows from theorem \ref{t.barriers-flat} that 
for every $a>0$, the pair of hypersurfaces
$\left(\Sigma_{a/(n-1)},\Sigma_a\right)$ is a pair of
$-1/a$-barriers. Hence,
theorem \ref{t.CMC-existence} and remark \ref{r.uniqueness} imply that the
hypersurface $\tau_{cmc}^{-1}(-1/a)$ is in the future of
$\Sigma_{a/(n-1)}=\tau^{-1}(a/(n-1))$ and in the past of 
$\Sigma_a=\tau^{-1}(a)$.
Equivalently, one has  
$$
\tau\leq -\frac{1}{\tau_{cmc}}\leq (n-1)\tau.
$$
The case where $(M,g)$ is future geodesically incomplete is similar
(except that $(M,g)$ is the quotient of a \emph{past complete} flat
regular domain $E^-(\Lambda)$, and that one has to consider the
\emph{reverse} cosmological time of $E^-(\Lambda)$). 

Finally, let us consider the case where $(M,g)$ is geodesically
complete. Then Theorem~\ref{teo.dscompact} states that up to a finite
covering $(M,g)$ is a
quotient of $\RR^{1,n-1}$ by a commutative subgroup $\Gamma$ of
$\mbox{Isom}(\RR^{1,n-1})$ generated by $n-1$ spacelike linearly independant
translations $t_{\overrightarrow{u_1}},\dots,t_{\overrightarrow{u_n}}$. 
Let $\overrightarrow{v}$ be any  (say future-directed) timelike vector. Then, 
for every $t\in\RR$, the affine plane 
$P_t:=t.\overrightarrow{v}+\RR.\overrightarrow{u_1}+\dots+\RR.\overrightarrow{u_n}$  
is $\Gamma$-invariant. Hence it induces a totally geodesic spacelike
hypersurface $\Sigma_t:=\Gamma\backslash P_t$ in
$M\simeq\Gamma\backslash\RR^{1,n-1}$. The family of hypersurfaces
$(\Sigma_t)_{t\in\RR}$ is a foliation of $M$ whose leaves are by totally geodesic 
(in particular, CMC) spacelike hypersurfaces. 

In order to complete the proof of Theorem~\ref{t.main-flat-case}, 
we only need to prove that in the case where  $(M,g)$ is geodesically complete, 
every CMC Cauchy hypersurface $\Sigma$ in $M$ is a leaf of the totally geodesic foliation 
$(\Sigma_t)_{t\in\RR}$ constructed above. 
Indeed, let  $t^-=\inf\{t\mbox{ such that }\Sigma\cap\Sigma_t\neq \emptyset\}$ and 
$t^+=\sup\{t\mbox{ such that }\Sigma\cap\Sigma_t\neq \emptyset\}$.
Then, $\Sigma$ is tangent to $\Sigma_{t^-}$ at some point and is in the future of 
$\Sigma_{t^-}$.
Hence, the maximum  principle (Proposition~\ref{p.Harnack}) implies that the
mean curvature of  
$\Sigma$ is smaller or equal than those of $\Sigma_{t^-}$, i.e. is
non-positive. Similarly,   
$\Sigma$ is tangent to $\Sigma_{t^+}$ at some point and is in the past of $\Sigma_{t^+}$, 
so by the maximum principle, the mean curvature of $\Sigma$ is non-negative. So, we know 
that the mean curvature of $\Sigma$ is equal to $0$. And now, we use the 
equality case of the maximum principle 
(see, e.g.,~\cite[Theorem 3.6]{gallomax}):
\textit{if $S$ and $S'$ are two CMC Cauchy hypersurfaces with the same mean curvature,  
which are tangent at some point, and such that $S'$ is in the future of $S$, then $S=S'$.} 
This shows that $\Sigma=\Sigma_{t^-}=\Sigma_{t^+}$; in particular, $\Sigma$ is a leaf of the 
totally geodesic foliation  $(\Sigma_t)_{t\in\RR}$.

\subsection*{Acknowledgements} The authors are grateful for the hospitality and support of the Isaac Newton
Institute in Cambridge, where part of the work on this paper was
performed. We thank Ralph Howard for some helpful remarks, and for
pointing out reference \cite{caffarelli:cabre}.

\providecommand{\bysame}{\leavevmode\hbox to3em{\hrulefill}\thinspace}

\bigskip

%


\begin{thebibliography}{99}



\bibitem{andflat} Lars Andersson, \textit{Constant mean curvature foliations of
flat space-times,\/} Commun. Anal. Geom., {\bf 10}
(2002), 1125--1150.



\bibitem{andersson:survey}
Lars Andersson, \emph{The global existence problem in general relativity}, The
  Einstein equations and the large scale behavior of gravitational fields,
  Birkh\"auser, Basel, 2004, pp.~71--120.


\bibitem{andersson:wedge}
\bysame, \emph{Constant mean curvature foliations of simplicial spacetimes},
  Comm. Anal. Geom. \textbf{13} (2005), 1--17.


\bibitem{ABBZ} L. Andersson, T. Barbot, F. B\'eguin and A. Zeghib, 
\textit{Cosmological time versus CMC time II: the deSitter and anti-deSitter
  cases}. In preparation.

\bibitem{ABBZ:asymptotic} L. Andersson, T. Barbot, F. B\'eguin and A. Zeghib, 
\textit{Asymptotic behaviour of CMC hypersurfaces in globally
hyperbolic flat spacetimes}. In preparation.

\bibitem{cosmic} 
L.~Andersson,~G.J.~Galloway, R.~Howard, \textit{The Cosmological Time
  Function,} Classical Quantum Gravity, \textbf{15}~(1998), 309--322.

  

\bibitem{gallomax} 
L.~Andersson,~G.J.~Galloway, R.~Howard, \textit{A strong maximum 
principle for weak solutions of quasi-linear elliptic equations 
with applications to Lorentzian and Riemannian geometry,\/}  
Comm. Pure Appl. Math., {\bf 51}  (1998),  no. 6, 581--624.

    
\bibitem{AndHow}
Lars Andersson and Ralph Howard, \emph{Comparison and rigidity theorems in
  semi-{R}iemannian geometry}, Comm. Anal. Geom. \textbf{6} (1998), no.~4,
  819--877.


\bibitem{AndMon}
L.  Andersson, V.  Moncrief, \textit{Elliptic-hyperbolic  systems and  the
Einstein  equations,} Ann.  Henri Poincar\'e,  \textbf{4} (2003),
no. 1, 1--34.  

    

\bibitem{AMT}
Lars Andersson, Vincent Moncrief, and Anthony~J. Tromba, \emph{On the global
  evolution problem in {$2+1$} gravity}, J. Geom. Phys. \textbf{23} (1997),
  no.~3-4, 191--205. 



\bibitem{barflat} 
T.~Barbot, \textit{Flat globally hyperbolic spacetimes,}
Journ. Geom. Phys., \textbf{53} (2005), 123--165.  



\bibitem{BBZ:2003}
Thierry Barbot, Fran{\c{c}}ois B{\'e}guin, and Abdelghani Zeghib,
  \emph{Feuilletages des espaces temps globalement hyperboliques par des
  hypersurfaces \`a courbure moyenne constante}, C. R. Math. Acad. Sci. Paris
  \textbf{336} (2003), no.~3, 245--250. 



\bibitem{BBZ} 
T.~Barbot, F.~B\'eguin, A.~Zeghib, \textit{Constant mean curvature
  foliations of globally hyperbolic spacetimes locally modelled on
  $AdS_3$,} math.MG/0412111, to appear in Geom. Ded. 



\bibitem{benedetti:guadagnini}
Riccardo Benedetti and Enore Guadagnini, \emph{Cosmological time in
  {$(2+1)$}-gravity}, Nuclear Phys. B \textbf{613} (2001), no.~1-2, 330--352.

 
\bibitem{benbon} 
R.~Benedetti,~F.~Bonsante, \textit{Canonical Wick rotations in
  3-dimensional gravity,} arXiv~: DG/0508485.  



\bibitem{bernal:sanchez:2003}
Antonio~N. Bernal and Miguel S{\'a}nchez, \emph{On smooth {C}auchy
  hypersurfaces and {G}eroch's splitting theorem}, Comm. Math. Phys.
  \textbf{243} (2003), no.~3, 461--470. 



\bibitem{bernal:sanchez:2005}
\bysame, \emph{Smoothness of time functions and the metric splitting of
  globally hyperbolic spacetimes}, Comm. Math. Phys. \textbf{257} (2005),
  no.~1, 43--50. 





\bibitem{Bons2} 
F.~Bonsante. \textit{Deforming the Minkowskian cone of a 
    closed hyperbolic manifold,\/} Ph. D. Thesis, Pisa, 2005. 



\bibitem{bonsante} 
F.~Bonsante, \textit{Flat spacetimes with compact hyperbolic Cauchy
    surface,\/} Journ. Diff. Geom.,~\textbf{69}(2005), 441--521.  



\bibitem{budic:etal}
Robert Budic, James Isenberg, Lee Lindblom, and Philip~B. Yasskin, \emph{On
  determination of {C}auchy surfaces from intrinsic properties}, Comm. Math.
  Phys. \textbf{61} (1978), no.~1, 87--95. 

\bibitem{caffarelli:cabre}
Luis~A. Caffarelli and Xavier Cabr{\'e}, \emph{Fully nonlinear elliptic
  equations}, American Mathematical Society Colloquium Publications, vol.~43,
  American Mathematical Society, Providence, RI, 1995. 

\bibitem{ChB:ruggeri}
Yvonne Choquet-Bruhat and Tommaso Ruggeri, \emph{Hyperbolicity of the $3+1$
  system of {E}instein equations}, Comm. Math. Phys. \textbf{89} (1983), no.~2,
  269--275.



\bibitem{chruscetal}
Piotr~T. Chrusciel, James Isenberg, and Daniel Pollack, \emph{Gluing initial
  data sets for general relativity}, Phys. Rev. Lett. \textbf{93} (2004),
  081101.




\bibitem{goldman:drumm}
Todd~A. Drumm and William~M. Goldman, \emph{The geometry of crooked planes},
  Topology \textbf{38} (1999), no.~2, 323--351.




\bibitem{Fischer-Moncrief}
Arthur~E. Fischer and Vincent Moncrief, \emph{Conformal volume collapse of
  3-manifolds and the reduced {E}instein flow}, Geometry, mechanics, and
  dynamics, Springer, New York, 2002, pp.~463--522. 



\bibitem{fried} 
D. Fried, \textit{Flat spacetimes,} J. Differential
  Geom.,~\textbf{26} (1987),  no. 3, 385--396. 



\bibitem{Ger1}
C. Gerhardt, \textit{$H$-surfaces in Lorentzian manifolds,}
Comm. Math. Phys.,  \textbf{89}  (1983), no. 4, 523--553.  



 

\bibitem{Ger3}
C. Gerhardt, \textit{On the CMC foliation of future ends of a spacetime.} 
arXiv:math.DG/0408197.




\bibitem{HawEll}
S.W. Hawking and G.F.R. Ellis, \textit{The large scale structure of
  space-time,} Cambridge University Press, 1973.



\bibitem{margulis}
G.~A. Margulis, \emph{Complete affine locally flat manifolds with a free
  fundamental group}, Zap. Nauchn. Sem. Leningrad. Otdel. Mat. Inst. Steklov.
  (LOMI) \textbf{134} (1984), 190--205, Automorphic functions and number
  theory, II.





\bibitem{mess} 
G.~Mess, \textit{Lorentz Spacetime of Constant Curvature,} IHES preprint,
1990.  A commented version of this preprint will appear in a special issue of Geom. Ded.



\bibitem{moncrief:teich}
Vincent Moncrief, \emph{Reduction of the {E}instein equations in {$2+1$}
  dimensions to a {H}amiltonian system over {T}eichm\"uller space}, J. Math.
  Phys. \textbf{30} (1989), no.~12, 2907--2914. 







\bibitem{rendall:survey}
Alan~D. Rendall, \emph{Theorems on existence and global dynamics for the
  {E}instein equations}, Living Rev. Relativ. \textbf{5} (2002), 2002--6, 62
  pp. (electronic).



 

\bibitem{Treibergs}
 A. Treibergs,  
 \textit{Entire Spacelike Hypersurfaces on Constant Mean Curvature in Minkowski Space,}
 Invent.  Math. \textbf{66}   (1985), 39-56.



\end{thebibliography}
\end{document}